\documentclass[12pt]{amsart}
\usepackage{amsmath, amssymb, amsfonts, amsthm, mathtools}
\usepackage{hyperref}
\usepackage{tabularx}
\usepackage{booktabs}
\usepackage{caption}
\usepackage{aurical}

\newtheorem{thm}{Theorem}[section]
\newtheorem{lem}[thm]{Lemma}
\newtheorem{cor}[thm]{Corollary}
\newtheorem{prop}[thm]{Proposition}
\newtheorem{ex}[thm]{Example}

\newtheorem{conp}[thm]{Pirashvili Conjecture}
\newtheorem*{prob*}{Open problem}

\theoremstyle{definition}
\newtheorem{defi}[thm]{Definition}

\theoremstyle{remark}
\newtheorem{rem}[thm]{Remark}
\newtheorem*{rem*}{Remark}

\DeclareMathOperator{\id}{id}

\DeclareMathOperator{\rad}{rad}
\DeclareMathOperator{\Hom}{Hom}

\newcommand{\kringel}{\mathbin{\raise0.5pt\hbox{$\scriptstyle\circ$}}}
\newcommand{\pkt}{\mathbin{\raise0.5pt\hbox{$\scriptstyle\bullet$}}}
\newcommand{\sq}{\mathbin{\raise0.5pt\hbox{$\scriptscriptstyle\square$}}}

\newcommand{\C}{\mathbb{C}}

\newcommand{\tr}{\mathop{\rm tr}}
\newcommand{\ad}{{\rm ad}}

\newcommand{\Der}{{\rm Der}}

\newcommand{\nil}{\mathop{\rm nil}}

\newcommand{\Lg}{\mathfrak{g}}

\newcommand{\Ll}{\mathfrak{l}}
\newcommand{\Ln}{\mathfrak{n}}

\newcommand{\Lr}{\mathfrak{r}}
\newcommand{\Ls}{\mathfrak{s}}

\newcommand{\om}{\omega}

\newcommand{\ra}{\rightarrow}

\renewcommand{\phi}{\varphi}

\begin{document}

% Ab hier duerfen Sie wieder.

\title[Sympathetic Lie algebras]{Sympathetic Lie algebras and adjoint cohomology for Lie algebras}

%  Die Kurzfassung kommt oben ueber die Seiten, sie steht in eckigen Klammern
%  Auch Autorennamen koennen eine Kurzfassung haben

\author[D. Burde]{Dietrich Burde}
\author[F. Wagemann]{Friedrich Wagemann}
\address{Fakult\"at f\"ur Mathematik\\
Universit\"at Wien\\
  Oskar-Morgenstern-Platz 1\\
  1090 Wien \\
  Austria}
\email{dietrich.burde@univie.ac.at}
\address{Laboratoire de math\'ematiques Jean Leray\\
  UMR 6629 du CNRS\\
  Universit\'e de Nantes \\
  2, rue de la Houssini\`ere, F-44322 Nantes Cedex 3 \\
  France}
\email{wagemann@math.univ-nantes.fr}

\date{\today}

\subjclass[2000]{Primary 17A32, Secondary 17B56}
\keywords{Lie algebra cohomology, sympathetic Lie algebras}

\begin{abstract}
We study sympathetic Lie algebras, namely perfect and complete Lie algebras. They arise among other things
in the study of adjoint Lie algebra cohomology. This is motivated by a conjecture of Pirashvili, which says that a
non-trivial finite-dimensional complex perfect Lie algebra is semisimple if and only if its adjoint cohomology vanishes.
We prove several results on sympathetic Lie algebras and the adjoint Lie algebra cohomology of Lie algebras in general, using
the Hochschild-Serre formula.
For certain semidirect products we obtain explicit results for the adjoint cohomology.
\end{abstract}

\maketitle

\section{Introduction}

% old
It is well-known that one can characterize finite-dimensional semisimple Lie algebras $\Lg$ 
over a field $K$ of characteristic zero by the vanishing of certain Lie algebra cohomology groups. 
For example, by Whitehead's first lemma, we have $H^1(\Lg,M)=0$ for every finite-dimensional $\Lg$-module $M$.
The converse statement is also true - any Lie algebra whose first cohomology with coefficients in any finite-dimensional 
module vanishes is semisimple. By Whitehead's second lemma, for a semisimple Lie algebra $\Lg$ we also have
$H^2(\Lg,M)=0$ for every finite-dimensional $\Lg$-module $M$. However, the converse is no longer true, see \cite{ZUS}.\\[0.2cm]
It has also been asked, whether or not the vanishing of the adjoint cohomology groups for $\Lg$ implies that $\Lg$ is semisimple.
This is not true in general, see for example the family of non-perfect reductive Lie algebras $\Lg$ given in
Example $\ref{2.8}$, which satisfies $H^n(\Lg,\Lg)=0$ for all $n\ge 0$. It is natural, however, to add the condition
$H^1(\Lg,\C)=0$ for the trivial module. Note that this is a strong condition on $\Lg$, which is equivalent to
$[\Lg,\Lg]=\Lg$, i.e., to $\Lg$ being perfect. \\[0.2cm]
The study of perfect Lie algebras with vanishing adjoint cohomology groups has already a long history.
In $1988$ Angelopoulos stated in \cite{A} that the question goes back to M. Flato some decades ago, who asked whether
semisimple Lie algebras  $\Lg$ are characterized by the vanishing conditions $H^1(\Lg,K)=H^1(\Lg,\Lg)=0$. Afterwards several
authors constructed complex non-semisimple Lie algebras $\Lg$ satisfying 
\[
H^1(\Lg,\C)=H^2(\Lg,\C)=H^0(\Lg,\Lg)=H^1(\Lg,\Lg)=H^2(\Lg,\Lg)=0,
\]
see \cite{A1,ABB,B1,B2}.
Benayadi \cite{B1} introduced {\em sympathetic} Lie algebras, i.e., Lie algebras which are perfect and complete, so that they
satisfy $H^1(\Lg,\C)=H^0(\Lg,\Lg)=H^1(\Lg,\Lg)=0$. He constructed in $1996$ a non-semisimple sympathetic Lie algebra $\Lg$ over the
complex numbers of dimension $25$. This Lie algebra has the lowest dimension among the known examples of complex non-semisimple
sympathetic Lie algebras.
The Lie algebra of Angelopoulos has dimension $35$, is sympathetic and satisfies $H^2(\Lg,\Lg)=0$. In this article we
will provide a basis and explicit Lie brackets for Benayadi's Lie algebra in dimension $25$, and show that
it satisfies $\dim H^2(\Lg,\Lg)=1$. \\[0.2cm]
In $2013$ T. Pirashvili \cite{P1} studied Lie algebra and Leibniz algebra cohomology and posed the conjecture,
that a non-trivial finite-dimensional complex Lie algebra $\Lg$ is semisimple if and only if it is perfect and
satisfies $H^n(\Lg,\Lg)=0$ for all $n\ge 0$. He called the conjecture the ``Weak Conjecture'', see page $1624$ in \cite{P1},
and also formulated a ``Strong Conjecture''. He proved one direction of the weak conjecture, namely that a semisimple
Lie algebra has vanishing adjoint cohomology and satisfies $H^1(\Lg,\C)=0$. \\[0.2cm]
The outline of this paper is as follows. In the second section we recall the definition and basic properties of sympathetic
Lie algebras and provide results on the adjoint cohomology of Lie algebras. We discuss the conjecture by Pirashvili
as stated above. We consider a few special cases and obtain some partial results. \\[0.2cm]
In the third section we study the adjoint cohomology of Lie algebras also for Lie algebras, which are not necessarily sympathetic.
Here we use the Hochschild-Serre formula and other tools from homological algebra. For Lie algebras $\Lg=\Ls\ltimes V$, where
$\Ls$ is semisimple, and $V$ is an $\Ls$-module we obtain non-vanishing results for $H^k(\Lg,\Lg)$. For $\Ls=\mathfrak{sl}_n(\C)$
and the natural representation $V$ of $\Ls$ we obtain an explicit result for all cohomology groups $H^k(\Lg,\Lg)$. \\[0.2cm]
In the fourth section we show that Benayadi's non-semisimple sympathetic Lie algebra $\Lg$ of 
dimension $25$ satisfies $\dim H^2(\Lg,\Lg)=1$. The crucial step here is to provide
explicit Lie brackets for $\Lg$ from the implicit construction in \cite{B2}. Then it is possible 
to compute the cohomology by using a computer algebra system. It follows that
this Lie algebra cannot be a counterexample to the Pirashvili conjecture.

\section{Sympathetic Lie algebras and a conjecture by Pirashvili}

We always assume that all Lie algebras are finite-dimensional, non-trivial and defined over the complex numbers.
Let us recall the notion of a sympathetic Lie algebra, see \cite{B2}.

\begin{defi}
A Lie algebra $\Lg$ is called {\em sympathetic}, if it is perfect and complete, i.e., if it satisfies $[\Lg,\Lg]=\Lg$ and
$Z(\Lg)=0$, $\Der(\Lg)=\ad(\Lg)$.
\end{defi}

Note that a sympathetic Lie algebra is {\em unimodular}, i.e., satisfies $\tr(\ad(x))=0$ for all $x$. \\[0.2cm]
We may characterize sympathetic Lie algebras in terms of Lie algebra cohomology.

\begin{lem}
A Lie algebra $\Lg$ is sympathetic if and only if $H^1(\Lg,\C)=H^0(\Lg,\Lg)=H^1(\Lg,\Lg)=0$.
\end{lem}  

\begin{proof}
We have $H^1(\Lg,\C)\cong (\Lg/[\Lg,\Lg])^{\ast}$ and hence a Lie algebra $\Lg$ is perfect if and only if $H^1(\Lg,\C)=0$.
Furthermore we have $H^0(\Lg,\Lg)\cong Z(\Lg)$ and $H^1(\Lg,\Lg)\cong \Der(\Lg)/\ad(\Lg)$, so that a Lie algebra $\Lg$ is complete
if and only if $H^0(\Lg,\Lg)=H^1(\Lg,\Lg)=0$.
\end{proof}

\begin{defi}
Let $\Lg$ be a Lie algebra. Denote by $\rad(\Lg)$ the solvable radical of $\Lg$ and by $\nil(\Lg)$ the nilradical.
\end{defi}  

The solvable radical of a sympathetic Lie algebra is nilpotent.

\begin{lem}\label{2.4}
Let $\Lg$ be a sympathetic Lie algebra. Then we have $\rad(\Lg)=\nil(\Lg)$.
\end{lem}  

\begin{proof}
Let us write $\Lr$ for $\rad(\Lg)$ and $\Ln$ for $\nil(\Lg)$, and let $\Lg=\Ls \ltimes \Lr$ be a Levi decomposition
of $\Lg$. Then we have a direct vector space sum $\Lg=\Ls\dotplus \Lr$.
Since $\Lg$ is perfect, we have 
\[
\Ls\dotplus \Lr = \Lg= [\Lg,\Lg] =[\Ls\dotplus\Lr,\Ls\dotplus \Lr] = \Ls\dotplus([\Ls,\Lr]+[\Lr,\Lr]).
\]
Because the sum is direct it follows that 
\[
\Lr =[\Ls,\Lr]+[\Lr,\Lr] \subseteq \ad (\Lg)(\Lr) \subseteq \Ln.
\]
The last inclusion follows from the fact that $D(\Lr)\subseteq \Ln$ holds for all derivations $D\in \Der(\Lg)$, hence
in particular for inner derivations $D=\ad(x)$. It follows that $\Lr=\Ln$ is nilpotent.
\end{proof}  

Sympathetic Lie algebras have been studied by many authors, see for example \cite{A,A1,ABB,B1,B2}. Several examples
of sympathetic non-semisimple Lie algebras were constructed. This is of particular interest
in connection with a conjecture by Pirashvili, the so-called {\em Weak Conjecture} from \cite{P1}.

\begin{conp}
A finite-dimensional complex Lie algebra is semisimple if and only if it satisfies $H^1(\Lg,\C)=0$ and $H^n(\Lg,\Lg)=0$ for all $n\ge 0$.
\end{conp}

\begin{rem}
The conjecture can also be formulated in terms of vanishing Leibniz homology with trivial conditions, i.e., that
$HL_n(\Lg,\C)=0$ for all $n\ge 1$. For the equivalence of these cohomological conditions see Lemma $4.2$ of \cite{P1}.
For results on cohomology of Leibniz algebras and Lie algebras see \cite{FW}.
\end{rem}  

Let us call these cohomological vanishing conditions the {\em Pirashvili conditions} for $\Lg$. It is known that
every semisimple Lie algebra satisfies these conditions, see \cite{P}. In fact, this follows from the 
first Whitehead Lemma and the following result of Carles \cite{C1}.

\begin{prop}\label{2.7}
Let $\Lg$ be a complete Lie algebra, whose solvable radical is abelian. Then we have $H^n(\Lg,\Lg)=0$ for all $n\ge 0$.
\end{prop}

The converse direction of Pirashvili's conjecture is still open. A Lie algebra satisfying the Pirashvili
conditions obviously is sympathetic, but we don't know, whether or not it is necessarily semisimple. \\[0.2cm]
It is also interesting to note that the conjecture need not be true if we omit the condition $H^1(\Lg,\C)=0$ from
the Pirashvili conditions.

\begin{ex}\label{2.8}
Let $\Lg=\mathfrak{aff}(\C^m)\cong \Lg\Ll_m(\C)\ltimes \C^m$ be the affine Lie algebra with $m\ge 1$.
Then we have $H^n(\Lg,\Lg)=0$ for all $n\ge 0$, but $\Lg$ is not semisimple.
\end{ex}

Indeed, it is easy to see that $\Lg$ is complete, i.e., that $H^0(\Lg,\Lg)=H^1(\Lg,\Lg)=0$, see 
Theorem $4.2$ in \cite{MEN}.
By Proposition $\ref{2.7}$ it follows that  $H^n(\Lg,\Lg)=0$ for all $n\ge 0$. \\
The following result shows that the Pirashvili conjecture is true for sympathetic Lie algebras, whose solvable radical is abelian.

\begin{lem}\label{2.9}
Let $\Lg$ be a sympathetic Lie algebra, whose solvable radical $\Lr$ is abelian. Then $\Lg$ is semisimple.
\end{lem}

\begin{proof}
Assume that $\Lr$ is abelian. Then $V=\Lr$ is an $m$-dimensional vector space. We have the Levi 
decomposition $\Lg=\Ls\ltimes V$ with a Levi subalgebra $\Ls$ and $V$ is an $\Ls$-module. Denote by $D$ the linear
map on $\Lg$ which is zero on $\Ls$ and the identity on $V$. We claim that it is a derivation of $\Lg$ with $\tr(D)=m$. 
The Lie bracket on $\Ls \dotplus V$ is given by
\[
[(x,v),(y,u)]=([x,y],x\cdot u-y\cdot v)
\]
for all $x,y\in \Ls$ and $u,v\in V$. Then $D([(x,v),(y,u)])=(0,x\cdot u-y\cdot v)$ and
\begin{align*}
[D((x,v)),(y,u)] & = [(0,v),(y,u)]=(0,-y\cdot v),\\
[(x,v),D((y,u))] & = [(x,v),(0,u)]=(0,x\cdot u).
\end{align*}
Since $\Lg$ is complete, it is an inner derivation. However, since $\Lg$ is perfect,
all adjoint operators $\ad (x)$ have zero trace. Hence $m=0$, so that $V=0$ and $\Lg$ is semisimple.
\end{proof}

It follows that a potential counterexample to Pirashvili's conjecture must have a nilpotent, non-abelian radical.

\begin{cor}
Let $\Lg$ be a non-semisimple Lie algebra, which satisfies the Pirashvili conditions. 
Then the solvable radical of $\Lg$ is nilpotent and non-abelian.
\end{cor}

\begin{proof}
The solvable radical is nilpotent by Lemma $\ref{2.4}$ and non-abelian by Lemma $\ref{2.9}$.
\end{proof}  

Similarly we can also obtain the following result.

\begin{prop}
Let $\Lg$ be a sympathetic Lie algebra with solvable radical $\Ln$. Suppose that we have a split short exact sequence
\[
0\ra Z(\Ln)\ra \Ln\ra \Ln/Z(\Ln)\ra 0.
\]
Then $\Lg$ is semisimple.
\end{prop}  
\begin{proof}
Let $\Lg\cong \Ls \ltimes \Ln$ be a Levi decomposition. The Lie bracket on the vector space $\Lg=\Ls\dotplus \Ln$ is given by
\[
[(s,n),(t,m)]=([s,t], [n,m]+ s\cdot m-t\cdot n)
\]
for $s,t\in \Ls$ and $n,m \in \Ln$. Since $\Ln$ is a central extension of $\Ln/Z(\Ln)$ by $Z(\Ln)$, the Lie bracket on the vector space
$\Ln=(\Ln/Z(\Ln))\dotplus Z(\Ln)$ is given by
\[
[n,m]=[(x,a),(y,b)]=([x,y],x.b-y.a+\om(x,y))
\]
for $x,y\in \Ln/Z(\Ln)$, $a,b\in Z(\Ln)$ and $\om \in H^2(\Ln/Z(\Ln),Z(\Ln))$.  Since the extension is central and split we may assume that $\om =0$. 
Writing $s\cdot m=s\cdot(x,a)=(s\cdot x,s\cdot a)$ and $t\cdot n=t\cdot (y,b)=(t\cdot y,t\cdot b)$, the Lie bracket on
$\Lg=\Ls\dotplus (\Ln/Z(\Ln))\dotplus Z(\Ln)$ becomes
\[
[(s,x,a),(t,y,b)]=([s,t],[x,y]+s\cdot x-t\cdot y,s.a-t.b)
\]
for $s,t\in \Ls$, $x,y\in \Ln/Z(\Ln)$ and $a,b\in Z(\Ln)$. \\[0.2cm]
Now define a linear map $D\colon \Lg\ra \Lg$ by $(s,x,a)\mapsto a$. It is a derivation of $\Lg$, because we have
\begin{align*}
D([(s,x,a),(t,y,b)]) & = (0,0,s.a-t.b),\\
[D((s,x,a),(t,y,b)] & = [(0,0,a),(t,y,b)] =(0,0,s.a),\\
[(s,x,a),D((t,y,b))] & = [(s,x,a),(0,0,b)] =(0,0,-t.b).
\end{align*}
Since $\Lg$ is complete, $D$ is an inner derivation. Let $D=\ad(z)$ and $m=\dim Z(\Ln)$. Now $\Lg$ is also complete so that
the adjoint operators have trace zero. Hence $\tr(D)=0$ and thus $m=0$ and $Z(\Ln)=0$. By Lemma $\ref{2.4}$, $\Ln$ is nilpotent.
Hence $Z(\Ln)=0$ implies that $\Ln=0$. This means that $\Lg$ is semisimple.
\end{proof}

\begin{cor}
Let $\Lg$ be a non-semisimple, sympathetic Lie algebra with solvable radical $\Ln$. Then $\Ln$ is nilpotent, non-abelian
and the extension $0\ra Z(\Ln)\ra \Ln\ra \Ln/Z(\Ln)\ra 0$ does not split.
\end{cor}

\section{Adjoint Lie algebra cohomology}

A well known construction for perfect but non-semisimple Lie algebras is the semidirect product $\Lg=\Ls\ltimes V$ of a semisimple
Lie algebra with a non-trivial simple $\Ls$-module $V$, considered as abelian Lie algebra, i.e., so that $\rad(\Lg)$
is abelian. Suppose that $\Lg$ is complete. Then $\Lg$ is sympathetic and hence semisimple by Lemma $\ref{2.9}$.
This is a contradiction. Thus $\Lg$ cannot be complete. In fact, this is true more generally, even if $\Lg$ is not perfect.

\begin{prop}\label{3.1}
Let $\Lg=\Ls\ltimes V$, where $\Ls$ is semisimple and $V$ is an $\Ls$-module. Then we have
\[
H^1(\Lg,\Lg)\cong H^1(V,V)^{\Ls}\cong\Hom_{\Ls}(V,V).
\]
In particular, $\dim H^1(\Lg,\Lg)\ge 1$ and $\Lg$ is not complete.
\end{prop}  

\begin{proof}
By Proposition $5.11$ in \cite{BU51} we have $H^1(\Lg,\Lg)\cong \Hom_{\Ls}(V,V)$. Since $V$ is an abelian Lie algebra and a trivial
$V$-module, $H^1(V,V)\cong \Hom(V,V)$, so that  $H^1(V,V)^{\Ls}\cong \Hom_{\Ls}(V,V)$.
Now we always have the identity in $\Hom_{\Ls}(V,V)$. Hence this space is at least $1$-dimensional.
\end{proof}  

\begin{cor}\label{3.2}
Let $\Lg=\Ls\ltimes V$, where $\Ls$ is semisimple and $V$ is a simple $\Ls$-module. Then we have
$\dim  H^1(\Lg,\Lg)=1$.
\end{cor}

\begin{proof}
By Schur's Lemma we have $\Hom_{\Ls}(V,V)\cong \C\cdot \id$. Hence the space is $1$-dimensional.
\end{proof}

\begin{prop}\label{3.3}
Let $\Lg=\Ls\ltimes V$, where $\Ls$ is semisimple and $V$ is an $\Ls$-module. Then we have an exact sequence
\[
0\ra H^1(V,V)^{\Ls}\ra H^1(V,\Lg)^{\Ls}\ra H^1(V,\Lg/V)^{\Ls}.
\]  
\end{prop}  

\begin{proof}
Consider the short exact sequence of $\Lg$-modules
\[
0\ra V\ra \Lg\ra \Lg/V\ra 0,
\]
which is also a short exact sequence of $V$-modules by restriction to $V\subset \Lg$. Here $V$ and $\Lg/V$ are trivial $V$-modules.
Applying the long exact sequence in cohomology we obtain
\[
\cdots \ra H^0(V,\Lg/V)\ra H^1(V,V)\ra H^1(V,\Lg)\ra H^1(V,\Lg/V)\ra \cdots
\]
Applying the functor of $\Ls$-invariants, which is exact on the subcategory of finite-dimensional $\Ls$-modules we obtain
\[
\cdots \ra H^0(V,\Lg/V)^{\Ls}\ra H^1(V,V)^{\Ls}\ra H^1(V,\Lg)^{\Ls}\ra H^1(V,\Lg/V)^{\Ls}\ra \cdots
\]
Since $H^0(V,\Lg/V)$ is the space of $V$-invariants of the trivial module $\Lg/V$, we obtain $H^0(V,\Lg/V)\cong \Lg/V$.
But we have $(\Lg/V)^{\Ls}=0$, because the quotient module $\Lg/V\cong \Ls$ does not contain non-zero $\Ls$-invariants.
Hence we have $H^0(V,\Lg/V)^{\Ls}=0$. This yields the claimed exact sequence.
\end{proof}

In particular the map $H^1(V,V)^{\Ls}\hookrightarrow  H^1(V,\Lg)^{\Ls}$ is injective, so that we have
\[
1\le \dim H^1(\Lg,\Lg)\le \dim H^1(V,\Lg)^{\Ls}.
\]

\begin{cor}
Let $\Lg=\Ls\ltimes V$, where $\Ls$ is semisimple and $V$ is an $\Ls$-module. Assume that $V$ does not contain any factor
isomorphic to a proper ideal of $\Ls$ in its decomposition of an $\Ls$-module. Then we have
\[
H^1(\Lg,\Lg) \cong H^1(V,\Lg)^{\Ls}.
\]
\end{cor}  

\begin{proof}
We have $H^1(V,\Lg/V)^{\Ls}\cong \Hom_{\Ls}(V,V/\Lg)$, because the Lie algebra $V$ is abelian and 
the $V$-module $\Lg/V$ is trivial. Both $\Lg/V\cong \Ls$ and $V$ decompose into direct factors, and 
by assumption they don't share an isomorphic factor. Hence we have $\Hom_{\Ls}(V,V/\Lg)=0$ 
and $H^1(V,\Lg/V)^{\Ls}=0$. So the exact sequence in Proposition $\ref{3.3}$ together with
Proposition $\ref{3.1}$ yields 
\[
H^1(V,\Lg)^{\Ls}\cong H^1(V,V)^{\Ls} \cong \Hom_{\Ls}(V,V)\cong H^1(\Lg,\Lg).
\]
\end{proof}  

The results on $H^1(\Lg,\Lg)$ can be generalized to higher cohomology groups as follows. Note that
we have for all $k$,
\begin{align*}
H^k(V,\Lg/V)^{\Ls} & \cong \Hom_{\Ls}(\Lambda^k(V),\Ls),\\
H^k(V,V)^{\Ls}  & \cong \Hom_{\Ls}(\Lambda^k(V),V).
\end{align*}  

\begin{prop}\label{3.5}
Let $\Lg=\Ls\ltimes V$, where $\Ls$ is semisimple and $V$ is an $\Ls$-module. Let $k\ge 1$ and suppose that the
$\Ls$-module $\Lambda^{k-1}(V)$ does not contain a submodule isomorphic to $\Ls$. Then we have an exact sequence
\[
0\ra H^k(V,V)^{\Ls}\ra H^k(V,\Lg)^{\Ls}\ra H^k(V,\Lg/V)^{\Ls}.
\]
Suppose that in addition the $\Ls$-module $\Lambda^k(V)$ does contain a submodule isomorphic to $V$. Then
we have $\dim H^k(\Lg,\Lg)\ge 1$.
\end{prop}  

\begin{proof}
As in the proof of Proposition $\ref{3.3}$ we have a long exact sequence
  \[
\cdots \ra H^{k-1}(V,\Lg/V)^{\Ls}\ra H^k(V,V)^{\Ls}\ra H^k(V,\Lg)^{\Ls}\ra H^k(V,\Lg/V)^{\Ls}\ra \cdots
\]
Here we have by the first assumption that
\[
 H^{k-1}(V,\Lg/V)^{\Ls}\cong \Hom_{\Ls}(\Lambda^{k-1}(V),\Ls)=0.
\]
So the first assertion follows. \\[0.2cm]
By the second assumption we have that $H^k(V,V)^{\Ls}\cong \Hom_{\Ls}(\Lambda^k(V),V,V)$
is nonzero. Thus the exact sequence implies that also $H^k(V,\Lg)^{\Ls}$ is nonzero.
Using Theorem $13$ of \cite{HS} on page $603$ for $M=\Lg$, $L=V$ and $K=\Ls$ we have
\[
H^k(\Lg,\Lg)\cong \bigoplus_{i+j=k}H^i(\Ls,\C)\otimes H^j(V,\Lg)^{\Lg}
\]
Note that we have $H^j(V,\Lg)^{\Lg}\cong H^j(V,\Lg)^{\Ls}$, see \cite{HS}, page 603, before Theorem $13$.
For $i=0$ and $j=k$ this direct sum contains the summand
\[
H^0(\Ls,\C)\otimes H^k(V,\Lg)^{\Ls} \cong \C\otimes H^k(V,\Lg)^{\Ls} \cong H^k(V,\Lg)^{\Ls},
\]
which is nonzero. Hence $H^k(\Lg,\Lg)$ is nonzero.
\end{proof}

In some cases we can explicitly compute all cohomology groups $H^k(\Lg,\Lg)$ by using the above arguments.

\begin{prop}\label{3.6}
Let $\Lg=\Ls\ltimes V$ with $\Ls=\mathfrak{sl}_n(\C)$, $n\ge 2$, and $V=L(\om_1)$ be the natural $\Ls$-module
of dimension $n$. Then we have for all $k\ge 1$
\[
  H^k(\Lg,\Lg)\cong H^{k-1}(\mathfrak{sl}_n(\C),\C).
\]
Here $H^{\ast}(\mathfrak{sl}_n(\C),\C)$ is isomorphic to the exterior algebra $\Lambda^*(c_3,c_5,\ldots,c_{2n-1})$ generated
by cocycles $c_{2i+1}$ for $i=1,\ldots ,n-1$.
\end{prop}  

\begin{proof}
The $\Ls$-module $\Lambda^k(V)$ has dimension $\binom{n}{k}$. It is irreducible for every $1\le k\le n$,
and it is different from the adjoint $\Ls$-module $\Ls$. Indeed, their dimensions are always different: 
suppose that $\binom{n}{k}=n^2-1$. It is well-known that
\[
\frac{n}{(n,k)}\mid \binom{n}{k}, \; 1\le k\le n,
\]
where $(n,k)=\gcd(n,k)$. If $k<n$, then $d=\frac{n}{(n,k)}$ is a divisor $d>1$ of $n$ and hence does not divide $n^2-1$.
Hence $\binom{n}{k}=n^2-1$ is impossible. For $k=n$ this is also impossible. Hence by Schur's Lemma we have
\[
H^k(V,\Lg/V)^{\Ls}\cong \Hom_{\Ls}(\Lambda^k(V),\Ls)=0
\]
for all $k\ge 1$. Then the long exact sequence from the proof of Proposition $\ref{3.5}$ yields
\[
H^k(V,\Lg)^{\Ls}\cong H^k(V,V)^{\Ls}\cong \Hom_{\Ls}(\Lambda^k(V),V)\cong \begin{cases} \C & \text{ for } k=1,n-1 \\
 0 & \text{ otherwise.} \end{cases}
\]
Thus the Hochschild-Serre formula yields
\begin{align*}
H^k(\Lg,\Lg) & \cong \bigoplus_{i+j=k}H^i(\Ls,\C)\otimes H^j(V,\Lg)^{\Ls} \\
             & \cong H^{k-1}(\Ls,\C) \otimes \C \\
             & \cong H^{k-1}(\mathfrak{sl}_n(\C),\C).
\end{align*}
It is well known that the cohomology $H^{\ell}(\mathfrak{sl}_n(\C),\C)$ is isomorphic to $\ell$-th component of
\[
\Lambda^*(c_3,c_5,\ldots,c_{2n-1})
\]
with generators $c_{2i+1}$. For a reference see \cite{ONV}, table $4$.
\end{proof}  

This yields, for example, with the $n$-dimensional natural $\mathfrak{sl}_n(\C)$-module $V(n)=L(\omega_1)$,
\[
H^{k+1}(\mathfrak{sl}_2(\C)\ltimes V(2),\mathfrak{sl}_2(\C)\ltimes V(2)) \cong H^k(\mathfrak{sl}_2(\C),\C)
=\begin{cases} \C, \text{ if } k=0,3 \\
0, \text{ otherwise} \end{cases}
\]  
and
\[
H^{k+1}(\mathfrak{sl}_3(\C)\ltimes V(3),\mathfrak{sl}_3(\C)\ltimes V(3)) \cong H^k(\mathfrak{sl}_3(\C),\C)
=\begin{cases} \C, \text{ if } k=0,3,5,8 \\
0, \text{ otherwise} \end{cases}
\]

\begin{rem}
One can also consider Proposition $3.6$ for $\Lg=\Ls\ltimes V$ with other classical Lie algebras $\Ls$ and their natural
representation $V$. We have used in the proof two facts depending on $\Ls$, namely that the exterior powers $\Lambda^k(V)$
are again irreducible, and that $\Lambda^k(V)$ and $\Ls$ are different as $\Ls$-modules for all $k\ge 1$.
Unfortunately this need not be true in general for simple Lie algebras $\Ls$ type $B_n,C_n,D_n$.
For example, for type $C_n$ the exterior powers $\Lambda^k(V)$ of dimension $\binom{2n}{k}$ are no longer irreducible
for $2\le k\le 2n-1$, see  \cite{FUH}, $\S 17.2$. And for types  $B_n$ and $D_n$, the $\Ls$-modules
$\Lambda^2(V)$ and $\Ls$ are no longer different. So one would need additional arguments for the computation of the
adjoint cohomology. 
\end{rem}

% Bemerkung von Friedrich
\begin{rem}
Using the Hochschild-Serre formula as in the proof of Proposition \ref{3.6}, one may also compute
the cohomology with trivial coefficients of a semi-direct product $\Lg=\Ls\ltimes V$ for a
general complex semisimple Lie algebra $\Ls$ and a general irreducible $\Ls$-module $V$.
In particular, if the $\Ls$-module $\Lambda^k(V)$ is irreducible for all $0\leq k\leq m$, for
$m:=\dim(V)$, then $H^k(V,\C)^\Ls=\Lambda^k(V)^\Ls=0$ for all $0<k<m$, and we obtain
$$H^n(\Lg,\C)=H^n(\Ls,\C)\oplus H^{n-m}(\Ls,\C).$$
This can also be used to compute the low degree {\em Leibniz cohomology} $HL^n(\Lg,\Lg)$
with adjoint coefficients in some cases. Consider the
$5$-dimensional Lie algebra $\Lg:=\mathfrak{sl}_2(\C)\ltimes V(2)$.
We obtain using the computational methods of Pirashvili in \cite{P} that
\[
\dim HL^n(\Lg,\Lg)=\begin{cases} 0, \text{ for } n=0,2, \\
1, \text{ for } n=1. \end{cases}
\]
Thus the Lie algebra $\Lg$ is rigid in the Leibniz sense.
Note that we used $H^\bullet(\Lg,\Lg)\cong H^\bullet(\Lg,\Lg^*)$ for the computation. This 
follows from the existence of an invariant, nondegenerate, symmetric bilinear form on $\Lg$.
\end{rem}

Finally we can also use the Hochschild-Serre formula to compute the adjoint cohomology $H^n(\Lg,\Lg)$ of certain semidirect products
$\Lg=\Ls\ltimes \Ln$ for the top degree $n$, i.e., for $n=\dim (\Lg)$.

\begin{prop}\label{3.8}
Let $\Lg$ be an $n$-dimensional Lie algebra with Levi decomposition $\Lg=\Ls\rtimes \rad(\Lg)$, where $\Ln=\rad(\Lg)$ is
nilpotent. Assume that the $\Ls$-module $\Ln/[\Ln,\Ln]$ does not contain the trivial $\Ls$-module $\C$. Then we have
\[
H^n(\Lg,\Lg)=0.
\]  
\end{prop}  

\begin{proof}
Let $\dim(\Ln)=m$ and $\dim(\Ls)=n-m$. By the Hochschild-Serre formula we have
\begin{align*}
H^n(\Lg,\Lg) & =\bigoplus_{p+q=n}H^p(\Ls,\C)\otimes H^q(\Ln,\Lg)^{\Ls} \\
             & = H^{n-m}(\Ls,\C)\otimes H^m(\Ln,\Lg)^{\Ls},
\end{align*}
because $\dim H^k(\Ln,\Lg)=0$ for all $k\ge m$, and $\dim H^k(\Ls,\C)=0$ for all $k\ge m-n$. To compute $H^m(\Ln,\Lg)^{\Ls}$ we 
use the long exact sequence as in the proof of Proposition $\ref{3.5}$, to obtain
\[
\cdots \ra H^m(\Ln,\Ln)^{\Ls}\ra H^m(\Ln,\Lg)^{\Ls}\ra H^m(\Ln,\Lg/\Ln)^{\Ls}\ra 0.
\]
Here $\Lg/\Ln$ is a trivial $\Ls$-module. Since $\Ln$ is nilpotent, we have $\dim H^m(\Ln,\C)\ge 1$ by Th\'eor\`eme $2$ of \cite{DIX}.
On the other hand, $\dim \Hom (\Lambda^m(\Ln),\C)= 1$, so that $\dim H^m(\Ln,\C)=1$. Therefore we have
\[
H^m(\Ln,\Lg/\Ln)\cong H^m(\Ln,\C)\otimes \Lg/\Ln \cong \Lg/\Ln\cong \Ls
\]
as $\Ls$-modules. Hence $H^m(\Ln,\Lg/\Ln)^{\Ls}\cong \Ls^{\Ls}=0$. \\[0.2cm]
To compute $H^m(\Ln,\Ln)^{\Ls}$, we use the Poincar\'e duality. Since $\Ln$ is unimodular, we obtain
\[
H^m(\Ln,\Ln)\cong H_0(\Ln,\Ln)\cong \Ln/[\Ln,\Ln].
\]
By assumption, the $\Ls$-module $\Ln/[\Ln,\Ln]$ does not contain the trivial module. Hence we have $H^m(\Ln,\Ln)^{\Ls}=0$.
Hence $H^m(\Ln,\Lg)^{\Ls}=0$, so that $H^n(\Lg,\Lg)=0$ by the Hochschild-Serre formula.
\end{proof}

\section{Cohomology of Benayadi's Lie algebra}

Benayadi constructed in \cite{B2} sympathetic non-semisimple Lie algebras of dimension $25$ by taking the vector space
\begin{align*}
\Lg & = \Ls\Ll_2(\C) \ltimes (V(7) \oplus V(5)\oplus V(7)\oplus V(3)),
\end{align*}
where $V(n)$ denotes the $n$-dimensional standard irreducible $\Ls\Ll_2(\C)$-module. He equipped $\Lg$ with a Lie bracket
such that the Lie brackets of $\mathfrak{sl}_2(\C)$ with $\Lg_1=V(7)$, $\Lg_2=V(5)$, $\Lg_3=V(7)$  and $\Lg_4=V(3)$ are given
by the action of $\mathfrak{sl}_2(\C)$ on $\Lg_i$, and such that
\[
[\Lg_1,\Lg_1]=\Lg_3,\; [\Lg_1,\Lg_2]=\Lg_4,\; [\Lg_2,\Lg_2]=\Lg_4,\; [\Lg_1,\Lg_3]=\Lg_4.
\]
We want to introduce a basis $\{e_1,\ldots ,e_{25}\}$ of $\Lg$, in order to obtain explicit Lie brackets.
Then the cohomology can be computed by a computer algebra system, e.g., with GAP. So fix a basis of $\Lg$, such that
$\{e_1,e_2,e_3\}$ is a basis of $\Ls\Ll_2(\C)$, $\{e_4,\ldots ,e_{10}\}$ is a basis of $\Lg_1$,
$\{e_{11},\ldots ,e_{15}\}$ is a basis of $\Lg_2$, $\{e_{16},\ldots ,e_{22}\}$ is a basis of $\Lg_3$ and
$\{e_{23},e_{24},e_{25}\}$ is a basis of $\Lg_4$. \\[0.2cm]
The action of $\Ls\Ll_2(\C)$ on $V(n)$ may be given by
\[
  \rho(e_1)=\begin{pmatrix} 0 & 1 & 0 & \cdots & 0 \\
    0 & 0 & 2 & \ddots & \vdots \\
  \vdots & \ddots & \ddots & \ddots & 0 \\
  \vdots & \ddots & \ddots & \ddots & n-1 \\
  0 & 0 & \cdots & 0 & 0    
\end{pmatrix},\quad
\rho(e_2)=\begin{pmatrix} 0 & 0 & \cdots & 0 & 0 \\
n-1 & 0 & \ddots & 0 & 0 \\
0 & \ddots & \ddots & \ddots & \vdots \\
\vdots & \ddots & 2 & 0 & 0 \\
0 & \cdots & 0 & 1 & 0
\end{pmatrix}
\]
and
\[
  \rho(e_3)=\begin{pmatrix} n-1 & 0 & \cdots  & 0 & 0 \\
  0 & n-3 & \cdots & 0 & 0 \\
  \vdots & \ddots & \ddots & \ddots & \vdots \\
  0 & 0 & \ddots & 3-n & 0 \\
  0 & 0 & \cdots & 0 & 1-n    
\end{pmatrix}
\]
So the nonzero brackets are determined as follows:\\[0.5cm]
$1.$ The brackets for $\Ls\Ll_2(\C):$
\begin{align*}
[e_1,e_2] & = e_3,      & [e_1,e_3] & =-2e_1,   & [e_2,e_3]& =2e_2.\\
\end{align*}
$2.$ The brackets between $\Ls\Ll_2(\C)$ and $\Lg_4:$
\begin{align*}
[e_1,e_{24}]  & = e_{23},   & [e_2,e_{23}] & =2e_{24}, & [e_3,e_{23}]& =2e_{23},\\
[e_1,e_{25}]  & = 2e_{24},  & [e_2,e_{24}] & =e_{25},  & [e_3,e_{25}]& =-2e_{25}, \\
\end{align*}
$3.$ The brackets between $\Ls\Ll_2(\C)$ and $\Lg_2:$
\begin{align*}
[e_1,e_{12}]  & = e_{11},   & [e_2,e_{11}] & =4e_{12}, & [e_3,e_{11}]& =4e_{11},\\
[e_1,e_{13}]  & = 2e_{12},  & [e_2,e_{12}] & =3e_{13}, & [e_3,e_{12}]& =2e_{12}, \\
[e_1,e_{14}]  & = 3e_{13},  & [e_2,e_{13}] & =2e_{14}, & [e_3,e_{14}]& =-2e_{14}, \\
[e_1,e_{15}]  & = 4e_{14},  & [e_2,e_{14}] & =e_{15},  & [e_3,e_{15}]& =-4e_{15}.\\
\end{align*}
$4.$ The brackets between $\Ls\Ll_2(\C)$ and $\Lg_3:$
\begin{align*}
[e_1,e_{17}]  & = e_{16},   & [e_2,e_{16}] & =6e_{17}, & [e_3,e_{16}]& =6e_{16},\\
[e_1,e_{18}]  & = 2e_{17},  & [e_2,e_{17}] & =5e_{18}, & [e_3,e_{17}]& =4e_{17}, \\
[e_1,e_{19}]  & = 3e_{18},  & [e_2,e_{18}] & =4e_{19}, & [e_3,e_{18}]& =2e_{18}, \\
[e_1,e_{20}]  & = 4e_{19},  & [e_2,e_{19}] & =3e_{20}, & [e_3,e_{20}]& =-2e_{20},\\
[e_1,e_{21}]  & = 5e_{20},  & [e_2,e_{20}] & =2e_{21}, & [e_3,e_{21}]& =-4e_{21}, \\
[e_1,e_{22}]  & = 6e_{21},  & [e_2,e_{21}] & =e_{22},  & [e_3,e_{22}]& =-6e_{22}.\\
\end{align*}
$5.$ The brackets between $\Ls\Ll_2(\C)$ and $\Lg_1:$
\begin{align*}
[e_1,e_{5}]  & = e_{4},   & [e_2,e_{4}] & =6e_{5}, & [e_3,e_{4}]& =6e_{4},\\
[e_1,e_{6}]  & = 2e_{5},  & [e_2,e_{5}] & =5e_{6}, & [e_3,e_{5}]& =4e_{5}, \\
[e_1,e_{7}]  & = 3e_{6},  & [e_2,e_{6}] & =4e_{7}, & [e_3,e_{6}]& =2e_{6}, \\
[e_1,e_{8}]  & = 4e_{7},  & [e_2,e_{7}] & =3e_{8}, & [e_3,e_{8}]& =-2e_{8},\\
[e_1,e_{9}]  & = 5e_{8},  & [e_2,e_{8}] & =2e_{9}, & [e_3,e_{9}]& =-4e_{9}, \\
[e_1,e_{10}]  & = 6e_{9},  & [e_2,e_{9}] & =e_{10},  & [e_3,e_{10}]& =-6e_{10}.\\
\end{align*}
$6.$ The brackets $[\Lg_1,\Lg_1]=\Lg_3:$
\begin{align*}
[e_4,e_5] & = a_1e_{16}+\cdots + a_7e_{22}, \\
[e_4,e_6] & = a_8e_{16}+\cdots + a_{14}e_{22}, \\
\cdots \hspace{0.3cm} & = \hspace{1.41cm} \cdots \\
[e_9,e_{10}] & = a_{141}e_{16}+\cdots + a_{147}e_{22}.\\
\end{align*}
$7.$ The brackets $[\Lg_1,\Lg_2]=\Lg_4:$
\begin{align*}
[e_4,e_{11}] & = b_1e_{23}+b_2e_{24}+b_3e_{25}, \\
[e_4,e_{12}] & = b_4e_{23}+b_5e_{24}+b_6e_{25}, \\
\cdots \hspace{0.3cm} & = \hspace{1.41cm} \cdots \\
[e_{10},e_{15}] & = b_{103}e_{23}+b_{104}e_{24}+b_{105}e_{25}.\\
\end{align*}
$8.$ The brackets $[\Lg_2,\Lg_2]=\Lg_4:$
\begin{align*}
[e_{11},e_{12}] & = c_1e_{23}+c_2e_{24}+c_3e_{25}, \\
[e_{11},e_{13}] & = c_4e_{23}+c_5e_{24}+c_6e_{25}, \\
\cdots \hspace{0.3cm} & = \hspace{1.41cm} \cdots \\
[e_{14},e_{15}] & = c_{28}e_{23}+c_{29}e_{24}+c_{30}e_{25}.\\
\end{align*}
$9.$ The brackets $[\Lg_1,\Lg_3]=\Lg_4:$
\begin{align*}
[e_4,e_{16}] & = d_1e_{23}+d_2e_{24}+d_3e_{25}, \\
[e_4,e_{17}] & = d_4e_{23}+d_5e_{24}+d_6e_{25}, \\
\cdots \hspace{0.3cm} & = \hspace{1.41cm} \cdots \\
[e_{10},e_{22}] & = d_{145}e_{23}+d_{146}e_{24}+d_{147}e_{25}.\\
\end{align*}

The Jacobi identity is equivalent to polynomial equations in the variables 
$a_i,b_i,c_i,d_i$. These can be easily solved by linear equations. It turns out that there are solutions.
The solution space only depends on the four nonzero parameters $a_{15},b_{13},c_7,d_{16}$.
We obtain a family of {\em Lie algebras $L(a,b,c,d)$} with
\[
(a_{15},b_{13},c_7,d_{16})=(3a,60b,2c,15d).
\]

The rewriting in terms of nonzero complex parameters $a,b,c,d$ is
only for our convenience, to avoid writing fractions.

\begin{prop}
The family of Lie algebras $L(a,b,c,d)$ has the following explicit Lie brackets with respect to the basis
$(e_1,\ldots ,e_{25})$.
\begin{align*}
[e_1,e_{2}]   & = e_{3},   & [e_2,e_{23}] & =2e_{24}, & [e_6,e_{8}]  & =-4ae_{19},\\
[e_1,e_{3}]   & = -2e_{1}, & [e_2,e_{24}] & =e_{25},  & [e_6,e_{10}] & =12ae_{21}, \\
[e_1,e_{5}]   & =  e_{4},  & [e_3,e_{4}]  & =6e_{4},  & [e_6,e_{13}] & =4be_{25}, \\
[e_1,e_{6}]   & = 2e_{5},  & [e_3,e_{5}]  & =4e_{5},  & [e_6,e_{14}] & =-8be_{24},\\
[e_1,e_{7}]   & = 3e_{6},  & [e_3,e_{6}]  & =2e_{6},  & [e_6,e_{15}] & =4be_{25},\\
[e_1,e_{8}]   & = 4e_{7},  & [e_3,e_{8}]  & =-2e_{8}, & [e_6,e_{19}] & =3de_{23}, \\
[e_1,e_{9}]   & = 5e_{8},  & [e_3,e_{9}]  & =-4e_{9}, & [e_6,e_{20}] & =2de_{24}, \\
[e_1,e_{10}]  & = 6e_{9},  & [e_3,e_{10}] & =-6e_{10},& [e_6,e_{21}] & =-5de_{25},
\end{align*}

\begin{align*}
[e_1,e_{12}]  & = e_{11},   & [e_3,e_{11}] & =4e_{11}, & [e_7,e_{8}] & =-3ae_{20},\\
[e_1,e_{13}]  & = 2e_{12},  & [e_3,e_{12}] & =2e_{12}, & [e_7,e_{9}] & =-3ae_{21}, \\
[e_1,e_{14}]  & = 3e_{13},  & [e_3,e_{14}] & =-2e_{14},& [e_7,e_{10}]& =3ae_{22}, \\
[e_1,e_{15}]  & = 4e_{14},  & [e_3,e_{15}] & =-4e_{15},& [e_7,e_{12}]& =-3be_{23},\\
[e_1,e_{17}]  & = e_{16},   & [e_3,e_{16}] & =6e_{16}, & [e_7,e_{13}]& =6be_{24},\\
[e_1,e_{18}]  & = 2e_{17},  & [e_3,e_{17}] & =4e_{17}, & [e_7,e_{14}]& =-3be_{25}, \\
[e_1,e_{19}]  & = 3e_{18},  & [e_3,e_{18}] & =2e_{18}, & [e_7,e_{18}]& =-3de_{23}, \\
[e_1,e_{20}]  & = 4e_{19},  & [e_3,e_{20}] & =-2e_{20},& [e_7,e_{20}]& =3de_{25},
\end{align*}

\begin{align*}
[e_1,e_{21}]  & = 5e_{20}, & [e_3,e_{21}] & =-4e_{21}, & [e_8,e_{9}] & =-2ae_{22},\\
[e_1,e_{22}]  & = 6e_{21}, & [e_3,e_{22}] & =-6e_{22}, & [e_8,e_{11}]& =4be_{23}, \\
[e_1,e_{24}]  & = e_{23},  & [e_3,e_{23}] & =2e_{23},  & [e_8,e_{12}]& =-8be_{24}, \\
[e_1,e_{25}]  & = 2e_{24}, & [e_3,e_{25}] & =-2e_{25}, & [e_8,e_{13}]& =4be_{25},\\
[e_2,e_{3}]   & = 2e_{2},  & [e_4,e_{7}]  & =3ae_{16}, & [e_8,e_{17}]& =5de_{23},\\
[e_2,e_{4}]   & = 6e_{5},  & [e_4,e_{8}]  & =12ae_{17},& [e_8,e_{18}]& =-2de_{24}, \\
[e_2,e_{5}]   & = 5e_{6},  & [e_4,e_{9}]  & =30ae_{18},& [e_8,e_{19}]& =-3de_{25}, \\
[e_2,e_{6}]   & = 4e_{7},  & [e_4,e_{10}] & =60ae_{19},& [e_9,e_{11}]& =20be_{24},
\end{align*}

\begin{align*}
[e_2,e_{7}]   & = 3e_{8},  & [e_4,e_{15}] & =60be_{23}, & [e_9,e_{12}]  & =-10be_{25},\\
[e_2,e_{8}]   & = 2e_{9},  & [e_4,e_{21}] & =15de_{23}, & [e_9,e_{16}]  & =-15de_{23}, \\
[e_2,e_{9}]   & = e_{10},  & [e_4,e_{22}] & =90de_{24}, & [e_9,e_{17}]  & =10de_{24}, \\
[e_2,e_{11}]  & = 4e_{12}, & [e_5,e_{6}]  & =-2ae_{16}, & [e_9,e_{18}]  & =5de_{25},\\
[e_2,e_{12}]  & = 3e_{13}, & [e_5,e_{7}]  & =-3ae_{17}, & [e_{10},e_{11}]& =60be_{25},\\
[e_2,e_{13}]  & = 2e_{14}, & [e_5,e_{9}]  & =10ae_{19}, & [e_{10},e_{16}]& =-90de_{24}, \\
[e_2,e_{14}]  & = e_{15},  & [e_5,e_{10}] & =30ae_{20}, & [e_{10},e_{17}]& =-15de_{25}, \\
[e_2,e_{16}]  & = 6e_{17}, & [e_5,e_{14}] & =-10be_{23},& [e_{11},e_{14}]& =2ce_{23},
\end{align*}

\begin{align*}
[e_2,e_{17}]  & = 5e_{18}, & [e_5,e_{15}] & =20be_{24}, & [e_{11},e_{15}]& =8ce_{24},\\
[e_2,e_{18}]  & = 4e_{19}, & [e_5,e_{20}] & =-5de_{23}, & [e_{12},e_{13}]& =-ce_{23}, \\
[e_2,e_{19}]  & = 3e_{20}, & [e_5,e_{21}] & =-10de_{24},& [e_{12},e_{14}]& =-ce_{24}, \\
[e_2,e_{20}]  & = 2e_{21}, & [e_5,e_{22}] & =15de_{25}, & [e_{12},e_{15}]& =2ce_{25},\\
[e_2,e_{21}]  & = e_{22},  & [e_6,e_{7}]  & =-3ae_{18}, & [e_{13},e_{14}]& =-ce_{25}. \\
\end{align*}
\end{prop}

It turns out that all Lie algebras $L(a,b,c,d)$ are isomorphic.

\begin{prop}
We have an isomorphism $L(a_1,b_1,c_1,d_1)\cong L(a_2,b_2,c_2,d_2)$ for all nonzero complex numbers
$a_1,a_2,b_1,b_2,c_1,c_2,d_1,d_2$.
\end{prop}

\begin{proof}
Let $\phi\colon L(a_1,b_1,c_1,d_1)\ra L(a_2,b_2,c_2,d_2)$ be the map given by $\phi(e_i)=\xi_i e_i$ for all $i$ with
$1\le i\le 25$. A direct computation shows that $\phi$ is a Lie algebra homomorphism if and only if
\begin{align*}
  \xi_1  = \xi_2 =   \xi_3 & = 1
\end{align*}

\begin{align*}
\xi_4  = \cdots =   \xi_{10} & = \frac{a_1b_2^2c_1d_1}{a_2b_1^2c_2d_2},\\[0.2cm]
\xi_{11}  = \cdots =  \xi_{15} & = \frac{a_1b_2^3c_1^2d_1}{a_2b_1^3c_2^2d_2}, \\[0.2cm]
\xi_{16}  = \cdots =  \xi_{22} & = \frac{a_1b_2^4c_1^2d_1^2}{a_2b_1^4c_2^2d_2^2}, \\[0.2cm]
\xi_{23}  = \xi_{24} =  \xi_{25} & = \frac{a_1^2b_2^6c_1^3d_1^2}{a_2^2b_1^6c_2^3d_2^2}. \\
\end{align*}
Obviously the determinant of the diagonal matrix associated to $\phi$ is nonzero. So the map is a Lie algebra isomorphism.
\end{proof}
Hence we may choose the parameters as  $(a_{15},b_{13},c_7, d_{16})=(3,60,2,15)$, namely by taking
\[
(a,b,c,d)=(1,1,1,1).
\]
Note that then all structure constants are integer valued. We call the Lie algebra
\[
L_{25}:=L(1,1,1,1)
\]
the {\em Benayadi Lie algebra}. It is uniquely determined up to isomorphism. Since we have explicit Lie brackets,
the cohomology can be easily computed by using a computer algebra system like GAP.

\begin{thm}
Benayadi's Lie algebra $\Lg=L_{25}$ satisfies $H^1(\Lg,\C)=H^0(\Lg,\Lg)=H^1(\Lg,\Lg)=0$ and $\dim H^2(\Lg,\Lg)=1$.
\end{thm}  

We also can determine some more adjoint cohomology of $\Lg$ without a computation. For example, applying Proposition
$\ref{3.8}$ gives the following result.

\begin{cor}
Benayadi's Lie algebra $\Lg=L_{25}$ satisfies $H^{25}(L,L)=0$.
\end{cor}  

\begin{proof}
For Benayadi's Lie algebra we have $\Ls=\mathfrak{sl}_2(\C)$ and $\Ln=\Lg_1\oplus \Lg_2\oplus \Lg_3\oplus \Lg_4$ with
\[
[\Lg_1,\Lg_1]=\Lg_3,\; [\Lg_1,\Lg_2]=\Lg_4,\; [\Lg_2,\Lg_2]=\Lg_4,\; [\Lg_1,\Lg_3]=\Lg_4.
\]
Hence we have $[\Ln,\Ln]\cong \Lg_3\oplus \Lg_4$, and $\Ln/[\Ln,\Ln]\cong \Lg_1\oplus \Lg_2$. Since this does not contain the trivial
$\Ls$-module $\C$, Proposition $\ref{3.8}$ implies that $H^{25}(L,L)=0$.
\end{proof}  

On the other hand, recall that $L$ is unimodular, since it is perfect. Hence we also obtain
\[
H^{25}(L,L)\cong H_0(L,L)\cong L/[L,L]=0
\]
directly by the Poincar\'e duality.

\section*{Acknowledgments}
Dietrich Burde is supported by the Austrian Science Foun\-da\-tion FWF, grant I 3248 and grant P 33811.

\end{document}